\documentclass{article}
\usepackage{amssymb,latexsym,amsmath,amsthm,mathrsfs}
\usepackage{graphicx}
\newcommand{\comment}[1]{}

\begin{document}
\title{On the infinity of infinities of orders of the infinitely large and infinitely
small\footnote{Originally published as
{\em De infinities infinitis gradibus tam infinite magnorum
quam infinite parvorum},
Acta Academiae scientiarum Imperialis Petropolitanae \textbf{2} (1780), no. I,
102--118.
E507 in the Enestr{\"o}m index.
Translated from the Latin by Jordan Bell,
Department of Mathematics, University of Toronto, Toronto, Canada.
Email: jordan.bell@gmail.com}}
\author{Leonhard Euler}
\date{}
\maketitle

1.\footnote{Translator: E507 is cited in pp. 84--86 of H. J. M. Bos,
{\em Differentials,
higher-order differentials and the derivative in the Leibnizian calculus}, Arch. Hist. Exact Sci.
\textbf{14} (1974), no. 1, 1--90. E507 is also mentioned on p. 700
of volume 4 of Moritz Cantor, {\em Vorlesungen \"uber Geschichte der
Mathematik}, 1894, and on p. 68 of Constantin Gutberlet, {\em Das unendliche:
Metaphysisch und mathematisch betrachtet}, G. Faber, 1878}
\footnote{Translator: Particularly useful references are G. H. Hardy, {\em Orders of infinity}, Cambridge Tracts in Mathematics, no. 12,
Cambridge University Press, 1910; p. 402 of Hardy's {\em A course of pure
mathematics}, seventh ed., Cambridge University Press, 1938; and p. 249 of Carl Boyer's {\em The history
of the calculus and its conceptual development}, Dover, 1959.}
If $x$ denotes an infinitely large quantity, then the geometric
progression
\[
1,x,xx,x^3,x^4,x^5,\, \textrm{etc.}
\]
is thus constituted that
each term is infinitely greater than the preceding, and indeed infinitely
smaller than the following. Whence if we think of
the power $x^{1000}$ 
as the last term of this progression, it will be possible to put 
a thousand orders of different infinite magnitudes between it and the
first term $1$; here we refer to
quantities having a finite ratio between themselves as being of the same order.
Yet still the number 1000 does not furnish all the intermediate orders between
$1$ and $x^{1000}$;
it should be noted that 
when we speak of the particular number 1000 here, any other number, however large,
can be put in its place.

2. As we just said, we are still far from
representing in this progression all the different intermediate orders 
between $1$ and $x^{1000}$.
For if we put\footnote{Translator: Original version has $z=y^{1000}$.}
$x=y^{1000}$ so that
\[
y=\sqrt[1000]{x},
\]
since $x$ is an infinite quantity even now $y$ will be an infinite
quantity; whence it follows that because again 
a thousand intermediate orders can be assigned
between $1$ and $y^{1000}$,
each of which is also infinitely greater than the preceding and
infinitely less than the following,
then even between unity and $x$ a thousand 
intermediate orders can be formed, even though before
$x$ had been the first order of infinity.
Indeed, 
in a similar way a thousand intermediate orders can again be assigned
between the preceding first order $x$ and the second order $xx$.
And it is thus between any two succeeding orders, which are all
such that each is infinitely greater than the preceding and infinitely
less than following.

3. But we still can't stop here. For since $y$ is an infinitely large
quantity, if we put
\[
y=z^{1000},
\]
then $z$ will still be an infinitely large
quantity; whence one sees that between 1 and $z^{1000}$, that is, between
1 and $y$, again a thousand orders of infinity can be constituted,
and we may proceed as far along in this fashion as we please, so that the number
of all the different orders can in fact be increased to infinity.

4. The same also holds for the infinitely small, in an inverted manner.
For if $x$ denotes an infinitely small quantity, then any term
of the geometric progression
\[
1,x,xx,x^3,\ldots,x^{1000}
\]
will be infinitely smaller
than the preceding but infinitely greater than the following, and hence
between 1 and $x^{1000}$ we obtain a thousand intermediate orders of
the infinitely small, all different;
for each is infinitely less than the preceding, and infinitely
greater than the following.

5. If now we again put $x=y^{1000}$, so that
\[
y=\sqrt[1000]{x},
\]
then $y$ will still be an infinitely small quantity; whence it is clear
that between 1 and $y^{1000}$, that is between 1 and $x$,
a thousand intermediate orders of the infinitely small can again be constituted,
which can also be done between $x$ and $xx$, and similarly between 
$xx$ and $x^3$, and so on in general between any two neighboring terms of
the preceding series. Since by further putting $y=z^{1000}$, $z$ is
even still an
infinitely small quantity, the number of different orders once again exceeds
one thousand, and this multiplication can be continued endlessly.

6. What has been deduced here from the consideration of powers is
indeed common knowledge, and can even be considered to be part of common Algebra;
but the higher Analysis also provides innumerable other 
orders of the infinitely large and infinitely small which cannot be put
into any of the orders that we have just related, no matter how many times
they are multiplied. Rather,  
they are found continually to be either infinitely greater or infinitely less than any of the preceding orders. Since I do not recall this being clearly
explained so far anywhere, this worthwhile task will be carefully
investigated here.\footnote{Translator: Euler also briefly discusses orders of infinity
in  \S 143, vol. 1 of his 1755
{\em Institutiones calculi differentialis} (E212).}

7. Now, these quantities occurring in higher Analysis can be put into
two classes, one of which is comprised of logarithms, and the other
exponential quantities. Let us first deal with
logarithms, and with $x$ denoting an infinitely large number it follows
that its logarithm will also be infinitely large. This is equally clear
whatever base we choose for logarithms, be it the common or hyperbolic one,
or any other type.\footnote{Translator: The hyperbolic logarithm is the
natural logarithm, the common logarithm is the logarithm with base 10.}

8. When $x$ is an infinitely large number, it is clear enough by itself
that its logarithm, that is,
\[
lx,
\]
will indeed be infinite, but still infinitely
less than the number $x$, and can be thus classed in a lower 
order. Since lower orders of $x$ can be represented by $x^{\frac{1}{n}}$,
namely with $n$ denoting some sufficiently large number, it
is hardly difficult to show that $lx$ is always infinitely less than
$x^{\frac{1}{n}}$,
however large a number is chosen for $n$.

9. In fact, it possible to demonstrate in the following way
that $x^{\frac{1}{n}}$ is
infinitely greater than $lx$ when $x=\infty$, or in other words
that the value of the fraction
\[
\frac{x^{\frac{1}{n}}}{lx}
\]
is infinitely
large.
For let us put this value $=v$, so that
\[
v=\frac{x^{\frac{1}{n}}}{lx},
\]
and set
\[
p=\frac{1}{lx} \quad \textrm{and} \quad q=\frac{1}{x^{\frac{1}{n}}},
\]
and so
\[
v=\frac{p}{q};
\]
both the numerator $p$ and the denominator $q$
of this fraction will be $=0$ in the case $x=\infty$.
Then, by the well known rule\footnote{Translator: Namely L'Hospital's rule. A footnote in the {\em Opera omnia} refers to Euler's {\em Institutiones calculi differentialis},
latter part of chapter XV, and also cites L'Hospital, {\em Analyse des infiniment petits}, Paris, 1696, p. 145} it will also be
\[
v=\frac{dp}{dq}.
\]
Therefore, since
\[
dp=-\frac{dx}{x(lx)^2} \quad \textrm{and} \quad
dq=-\frac{dx}{nx^{\frac{1}{n}+1}},
\]
we will have
\[
v=\frac{nx^{\frac{1}{n}}}{(lx)^2},
\]
which should be equal to the preceding value $\frac{x^{\frac{1}{n}}}{lx}$.
Truly then, taking the square of this will give
\[
vv=\frac{x^{\frac{2}{n}}}{(lx)^2},
\]
and dividing this value by the previous yields\footnote{Translator: In fact
$\frac{x^{\frac{2}{n}}}{(lx)^2}$ divided by $\frac{nx^{\frac{1}{n}}}{(lx)^2}$
is equal to $\frac{x^{\frac{1}{n}}}{n}$. But $n$ is a fixed positive integer, so $v$ is
still infinite.}
\[
v=nx^{\frac{1}{n}}.
\]
Since this is clearly infinite, it is thus apparent that 
$\frac{x^{\frac{1}{n}}}{lx}$ is an infinitely large quantity,
that is, that $lx$ is infinitely smaller than $x^{\frac{1}{n}}$ no matter
how large a number is taken for $n$.

10. It is therefore clear that
if $x=\infty$, its logarithm $lx$ cannot be grouped with any of the above
orders of infinity,
no matter how narrowly we divide these orders.
Because of this, a new
order of infinity, appropriate for the classification of the logarithm
$lx$,
must be constituted here, to which the power $x^{\frac{1}{n}}$
approaches continually closer the greater the number that is put
for $n$.
Yet the case where $n=\infty$ does not thereby satisfy this,
for when $\frac{1}{n}=0$ it makes $x^{\frac{1}{n}}=1$,\footnote{Translator:
As $x \to +0, x^x \to 0$.}
while $lx$ will be infinite; indeed it should
be noted that the demonstration 
given above produced $nx^{\frac{1}{n}}$,
which by taking $n=\infty$ nonetheless  
yields $v=n$, which is thus infinite.\footnote{Translator: As we noted earlier,
$v=\frac{x^{\frac{1}{n}}}{n}$ not $nx^{\frac{1}{n}}$. In fact,
if $x^{\frac{1}{n}}=1$ then $v=0$, not infinity.}

11. Therefore, since $lx$ constitutes as it were
the lowest order of all
the infinitely large quantities, it is evident that 
the number of orders that we established above,
which was already seen to be infinite,
can be further increased infinitely. For if one considers the order
designated by the power $x^\alpha$, then the formula
\[
x^\alpha lx
\]
will be infinitely greater than $x^\alpha$; yet 
when the exponent $\alpha$ is increased by a fraction $\frac{1}{n}$,
then the formula $x^\alpha lx$ will certainly be infinitely smaller
than $x^{\alpha+\frac{1}{n}}$, and so an order must be constituted between $x^\alpha$
and $x^{\alpha+\frac{1}{n}}$.

12. Yet the multitude of all the different orders is by no means exhausted
in this way. For $(lx)^2$ will be infinitely greater than $lx$ and thus
ought to form its own order, which will still be less than the powers
$x^{\frac{1}{n}}$ however large the number $n$ is taken to be.
All the different powers of $lx$ yield their own kind of infinity which has
to be extended even to fractional powers. For, since $(lx)^{\frac{\alpha}{\beta}}$
is certainly infinitely greater than $(lx)^{\frac{\alpha}{\beta}-\frac{1}{n}}$
yet infinitely less than $(lx)^{\frac{\alpha}{\beta}+\frac{1}{n}}$,
it should form its own order.
As many new cases will arise if in addition we multiply by a power of
$x$. Namely, the formula $x^a(lx)^{\frac{\alpha}{\beta}}$
is infinitely greater than $x^a(lx)^{\frac{\alpha}{\beta}-\frac{1}{n}}$,
but infinitely less than
$x^a(lx)^{\frac{\alpha}{\beta}+\frac{1}{n}}$.

13. Yet still in this manner we have not enumerated all the orders of
the infinite.
For since $lx$ is an infinitely large quantity, its own logarithm $llx$ will
still be infinite though infinitely less than $lx$;
whence one sees that again infinitely many new orders of the infinite
have to be established from the formula
\[
llx \quad \textrm{and its powers} \quad (llx)^{\frac{\alpha}{\beta}},
\]
first if this formula is combined not only with powers of $lx$ but also
with powers of $x$; and the same consideration can be extended further
to formulas
\[
lllx, \quad llllx, \quad \textrm{etc.}
\]

14. This immense multitude of orders also occurs among the infinitely small,
since of course they can be seen as reciprocals of infinite magnitudes.
For if unity is divided by any infinity $\infty$, namely $\frac{1}{\infty}$,
it should be considered to constitute a particular order of the infinitely
small.
Thus if $x$ is an infinite quantity, not only will the series
\[
\frac{1}{x}, \quad \frac{1}{xx}, \quad \frac{1}{x^3}, \quad \frac{1}{x^4}, \quad \textrm{etc.}
\]
provide
infinitely many orders of the infinitely small, but also
the series
\[
\frac{1}{lx}, \quad \frac{1}{(lx)^2}, \quad \frac{1}{(lx)^3},
\quad \frac{1}{(lx)^4}, \quad  \textrm{etc.}
\]
together with all the powers of all the terms will yield new orders
of the infinitely small; then also, the series
\[
\frac{1}{llx},\quad
\frac{1}{(llx)^2},\quad \frac{1}{(llx)^3},\quad \textrm{etc.}
\]
and also all the further ones
where more logarithms are taken such as $lllx, llllx$, etc.,
and thus this
multitude is augmented infinitely.

15. What has been propounded so far for logarithms can similarly be extended
to exponential quantities, from which innumerable new orders of both the infinitely
large and the infinitely small can be likewise constituted that will be
entirely distinct from the preceding. For if, as has been the case so far,
$x$ denotes an infinite number, it is well known that the value of the power
$a^x$ will also be infinitely large exactly when $a$ is a number greater than
unity, while if $a<1$ then the power $a^x$ will exhibit an infinitely
small quantity. Now, let us first consider the infinitely large, by taking
$a>1$, 
and it's clear that the power $a^x$ is not only
infinitely greater than its exponent, but it can even be demonstrated
that the quantity $a^x$ is infinitely greater than the power $x^n$
however large the exponent $n$ may be. The demonstration can be obtained
in the following way.

16. Let us put
\[
\frac{a^x}{x^n}=v,
\]
and let
\[
p=\frac{1}{x^n} \quad \textrm{and} \quad q=\frac{1}{a^x},
\]
so we have $v=\frac{p}{q}$. Both the numerator
$p$ and the denominator $q$ of this fraction vanish in the case $x=\infty$,
and so
\[
v=\frac{dp}{dq}.
\]
Now,
\[
dp=-\frac{ndx}{x^{n+1}} \quad \textrm{and} \quad dq=-\frac{dxla}{a^x},
\]
whence
\[
v=\frac{na^x}{x^{n+1}la}.
\]
Indeed this formula is much more complicated than the given
$v=\frac{a^x}{x^n}$, so it might seem that nothing could be concluded.
But from the comparison of these formulas the true value of $v$
will be concluded. For, from the first
\[
v^{n+1}=\frac{a^{x(n+1)}}{x^{n(n+1)}}
\]
and from the latter
\[
v^n=\frac{n^n}{a^{nx}}{x^{n(n+1)}(la)^n},
\]
and dividing the first value by the latter gives 
\[
v=\frac{a^x(la)^n}{n^n},
\]
which is clearly an infinite value.
Thus it has been demonstrated that the formula $a^x$ is always infinitely
greater than $x^n$, no matter how large the exponent $n$ is taken, providing
$a>1$.
It is apparent then that the exponential quantity $a^x$ is infinitely
greater than all the orders of the infinite arising from the power
$x^n$.
Hence, even though $x$ is an infinite quantity, still
all the fractions
\[
\frac{a^x}{x}, \quad \frac{a^x}{xx}, \quad
\frac{a^x}{x^3}, \quad \textrm{and in general} \quad \frac{a^x}{x^n}
\]
forms
an order of infinity that exceeds all orders of infinity from the first
class. It is clear that this also happens for the formulas
$a^{\alpha x}$,
providing $\alpha>0$, and it will even happen for the formulas $a^{\alpha x^\beta}$, as long as positive values are taken for the letters $\alpha$
and $\beta$; these are therefore infinitely many infinities that are higher
than the powers of $x$, however large those may be.

17. Besides this, it should also be noted that even if
the formula $a^x$ belongs to an infinitely high order of the infinite,
still, however little the value of the letter $\alpha$ is increased,
the value of this formula will rise to a higher infinity. For if
$b>a$ then the formula $a^x$ will be to the formula $b^x$ as $1$ is to $(\frac{b}{a})^x$, that is, as $1$ to infinity, an order of infinitesimals.

18. Whenever $a>1$, then all the powers of $a^x$ can be translated into
powers of the fixed number $e$, whose hyperbolic logarithm is $=1$,
because
\[
a^x=e^{xla},
\]
and thus all infinities of this type can be represented in the form
$e^{\alpha x^\beta}$,
supposing $\alpha>0$ and $\beta>0$.
Then the formula
\[
\frac{e^{\alpha x^\beta}}{x^n}
\] 
will belong to an infinitely lower order of the infinite, but on the other hand
this case be more than made up for if we write $e^{\alpha x^\beta}$ in place
of $\alpha x^\beta$, which leads to this form
\[
e^{e^{\alpha x^\beta}}
\]
and we can repeat this process arbitrarily many times.

19. This can all transferred in an inverted manner to the infinitely small,
which we shall now consider in some detail. Thus let $x$ denote an infinitely
small quantity, whose powers
\[
x^\alpha
\]
therefore form innumerable
orders of the infinitely small; if the exponent $\alpha$ is increased by however small an amount they become infinitely smaller. However all of these
orders can be included in the first class of the infinitely small, if
the exponents $\alpha$ are understood to take all positive values.

20. Let us refer to as the second class the infinitely small which arise from
logarithms. For whenever $l\frac{1}{x}$ is infinite, it's reciprocal
\[
\frac{1}{l\frac{1}{x}}
\]
will be infinitely small. Now, for convenience let us put
\[
l\frac{1}{x}=u,
\]
so that this form becomes $\frac{1}{u}$, which will be
infinitely small such that it is infinitely greater than all the infinitely small
in the first class. The formulas
\[
\frac{1}{uu}, \quad \frac{1}{u^3}, \quad \frac{1}{u^4}, \quad \textrm{etc. and in general} \quad \frac{1}{u^\alpha}
\]
also belong to this class.
The forms
\[
\frac{x^\alpha}{u^\beta}
\]
are the next to be considered,
which are combinations from the first and second classes.
Also indeed, since $lu$ is infinitely large but infinitely less than $u$, its
reciprocal
\[
\frac{1}{lu}
\]
will be infinitely small but infinitely greater than $\frac{1}{u}$. Likewise,
the formulas
\[
\frac{1}{llu} \quad \textrm{and} \quad \frac{1}{lllu}
\]
will be continually
infinitely larger than the preceding. By combining them with the previous,
innumerable new orders of the infinitely small can be constituted which one
can by no means enumerate.

21. In particular it ought to be noted that even though $u=l\frac{1}{x}$ is infinitely
large, yet all the products $x^nu$ will be infinitely small, when $n>0$.
And even though this follows from the preceding, it can be succinctly demonstrated thus.
Put
\[
x^nu=v,
\]
and let
\[
x^n=p \quad \textrm{and} \quad \frac{1}{u}=q,
\]
so that
\[
v=\frac{p}{q}.
\]
Both the numerator and the denominator of this fraction vanish when $x=0$,
whence we further have
\[
v=\frac{dp}{dq}.
\]
Indeed
\[
dp=nx^{n-1}dx,
\]
and
because $u=l\frac{1}{x}$, i.e. $u=-lx$, it will be
\[
du=-\frac{dx}{x},
\]
and hence
\[
dq=\frac{dx}{xuu},
\]
from which we get $v=nx^nuu$.
From the first value we have $vv=x^{2n}uu$,
which divided by that just found gives
\[
v=\frac{x^n}{n},
\]
from which it is clear
that the value of $v$ is infinitely small, which also holds for the formula
\[
x^n u^\alpha,
\]
and hence not only when $\alpha$ is a positive number but also when it is negative,
since the formula $\frac{x^n}{u^m}$ is by itself infinitely
small.\footnote{Translator: Since $x$ is infinitely  small and
$u$ is infinitely large,
$\frac{x^n}{u^m}$ is certainly
infinitely small. I don't see why Euler writes $u^m$ instead of $u^\alpha$.}

22. Besides these two classes of the infinitely small, the exponential quantities
offer a third class to us. For when $x=0$, the formula $e^{\frac{1}{x}}$
exhibits as it were an infinite magnitude of the highest degree,
whose reciprocal
\[
\frac{1}{e^{\frac{1}{x}}}=e^{\frac{-1}{x}}
\]
expresses the infinitely small of the highest degree, which namely
will be infinitely smaller than any of the infinitely small from
the first class, and in fact this should hold for the general form
\[
\frac{1}{e^{\frac{\alpha}{x^\beta}}}.
\]
Now, for the sake of brevity let us put
\[
e^{\frac{\alpha}{x^\beta}}=v,
\]
so that we can express the infinitely small counterpart in the form
$\frac{1}{v}$.
Then since $lv=\frac{\alpha}{x^\beta}$, by differentiating we'll have
\[
\frac{dv}{v}=-\frac{\alpha \beta dx}{x^{\beta+1}}, \quad \textrm{and hence}
\quad
dv=-\frac{\alpha\beta vdx}{x^{\beta+1}}.
\]
In should be noted here in particular that even if $x$ is a vanishing quantity,
still the formulas $\frac{1}{x^nv}$ will express the infinitely small
of the highest degree.\footnote{Translator: Perhaps because
by L'Hospital's rule, $\lim_{x \to 0} \frac{1}{x^nv}=\lim_{x \to 0}
\frac{n}{\alpha\beta}x^{\beta-n}\cdot \frac{1}{v}$.
Doing this $k$ times will yield
$\frac{n(n-\beta)(n-2\beta)\cdots (n-(k-1)\beta)}{(\alpha\beta)^{k+1}}
x^{n-k\beta}\cdot \frac{1}{v}$, 
which for $k=\lceil \frac{n}{\beta} \rceil$ has the limit $0$ as $x\to 0$.}

23. With these classes now constituted, useful aids both for differentiating
and for integrating such infinitely small quantities can be found. For instance,
if for the first class we put
\[
ax^\alpha=y,
\]
it will be
\[
\frac{dy}{dx}\alpha ax^{\alpha-1}
\]
and
\[
\int ydx=\frac{a}{\alpha+1}x^{\alpha+1},
\]
and it is clear that this integral will be infinitely less than $y$, providing
that infinity is greater than the differential $\frac{dy}{dx}$; indeed this
can be made infinitely large if $\alpha<1$ and should for this purpose
be thought of as
belonging to the infinitely small of other classes.

24. Now let's consider the infinitely small of the second class,
and put
\[
l\frac{1}{x}=u,
\]
so that
\[
du=-\frac{dx}{x}.
\]
Let us set
\[
y=ax^\alpha u^m,
\]
where $\alpha>0$, 
$m$ indeed either a positive or negative number since in either case
this formula is infinitely small. Then it will happen that
\[
\frac{dy}{dx}=\alpha ax^{\alpha-1}u^m-amx^{\alpha-1}u^{m-1}
=ax^{\alpha-1}u^{m-1}(\alpha u-m),
\]
and because $u$ is infinite, by neglecting the latter term $-m$ it will follow that
\[
\frac{dy}{dx}=\alpha ax^{\alpha-1}u^m.
\]
Multiplying by
$dx$ and integrating gives
\[
\int a\alpha x^{\alpha-1}u^m=dx=y=ax^\alpha u^m,
\]
from which arises this rather memorable integration:
\[
\int x^{\alpha-1}u^m dx=\frac{1}{\alpha} x^\alpha u^m,
\]
or, writing $\beta$ in place of $\alpha-1$, it will be
\[
\int x^\beta u^m dx=\frac{1}{\beta+1}x^{\beta+1}u^m.
\]

25. Then if take the curved line whose ordinate
\[
y=ax^\beta u^m
\]
corresponds
to the abscissa $x$, where we let $\beta>1$ and the exponent $m$ is either
positive or negative,
such that the initial ordinate of the curve, where $x=0$, vanishes,
the area of this curve corresponding to the infinitely small abscissa $x$
will be
\[
\int ydx=\frac{a}{\beta+1}x^{\beta+1}u^m=\frac{1}{\beta+1}xy;
\]
namely it will be equal to a rectangle formed between the abscissa
$x$ and the ordinate $y$, divided by $\beta+1$,
which is all the more remarkable as the formula $x^\beta u^m dx$ cannot
be integrated except for that small portion of cases in which the exponent
$m$ is a positive integer.

26. Now let's consider likewise 
the infinitely small of the third class, and for the sake of brevity let us
put as above
\[
e^{\frac{\alpha}{x^\beta}}=v,
\]
so that it is
\[
dv=-\frac{\alpha\beta vdx}{x^{\beta+1}},
\]
and, as we have seen, the formula
\[
\frac{x^m}{v}
\]
will always be an infinitely
small quantity whether the exponent $m$ is positive or negative.\footnote{Translator: See \S 22.}
And if one puts
\[
\frac{x^m}{v}=z,
\]
it will be
\[
\frac{dz}{dx}=\frac{mx^{m-1}+\alpha\beta x^{m-\beta-1}}{v}=\frac{x^{m-\beta-1}}{v}(mx^\beta+\alpha\beta)
\]
where, because $mx^\beta$ vanishes before $\alpha\beta$, it will be
\[
\frac{dz}{dx}=\frac{\alpha\beta x^{m-\beta-1}}{v}
\]
whence by integrating in turn we'll have
\[
z=\alpha\beta \int \frac{x^{m-\beta-1}dx}{v}=\frac{x^m}{v},
\]
and then if we write $n$ in place of $m-\beta-1$, so $n$ can thus be any positive or negative number, it will always be
\[
\int \frac{x^ndx}{v}=\frac{1}{\alpha\beta} \frac{x^{n+\beta+1}}{v};
\]
this integration is correct as long as $x$ is infinitely small, though
this differential expression altogether refuses integration.

27. Now if the curved line is taken whose ordinate corresponding
to the abscissa $x$ is
\[
y=\frac{ax^n}{v},
\]
with
\[
v=e^{\frac{\alpha}{x^\beta}},
\]
where $\alpha$ and $\beta$ are positive numbers, and indeed the exponent $n$
can be either positive or negative, the ordinate of this curve
will vanish at the beginning, where $x=0$. The area of this curve corresponding
to an infinitely small abscissa $x$ will be
\[
\int ydx=\frac{a}{\alpha\beta} \frac{x^{n+\beta+1}}{v}=\frac{1}{\alpha\beta}
x^{\beta+1}y,
\] 
and thus if
\[
y=\frac{ax^n}{e^{\frac{1}{x}}},
\]
where $\alpha=1$ and $\beta=1$,
it will be
\[
\int ydx=xxy.
\]
That is: the area of the curve will be
equal with the rectangle formed from the square of the abscissa, and the ordinate.

28. Now, if we then look for a curve whose area in general will be
\[
\int y dx=xxy,
\]
this leads to the differential equation
\[
ydx=2xydx+xxdy,
\] 
whence we get
\[
\frac{dy}{y}=\frac{dx(1-2x}{xx}
\]
and by then integrating
\[
ly=-\frac{1}{x}-2lx,
\]
and writing this without logarithms,
\[
y=\frac{a}{xxe^{\frac{1}{x}}},
\]
which is contained in the given form if one takes $n=-2$,
and indeed the above integration will work whenever $x$ is infinitely small.

29. The last integration will even work if the infinitely small
quantity also involves any member of the second class. For let
\[
l\frac{1}{x}=u,
\]
and let us put
\[
z=\frac{ax^mu^n}{v}
\]
(where the exponents
$m$ and $n$ can be taken to be either positive or negative, since these
quantities will always be infinitely small when $v=e^{\frac{\alpha}{x^\beta}}$).
So that we can more easily unravel the value $\frac{dz}{dx}$, let us take
logarithms, and it will be
\[
lz=la+mlx+nlu-lv
\]
and hence
\[
\frac{dz}{z}=\frac{mdx}{x}+\frac{ndu}{u}-\frac{dv}{v}.
\]
Indeed because
\[
du=-\frac{dx}{x} \quad \textrm{and} \quad dv=-\frac{\alpha\beta dx}{x^{\beta+1}}v,
\]
substituting these two values into the above expression, it will
take on the following form:
\[
\frac{dz}{zdx}=\frac{m}{x}-\frac{n}{ux}+\frac{\alpha\beta}{x^{\beta+1}},
\]
where since $\beta+1>1$, both the prior terms vanish before the third, and
thus it will neatly be
\[
\frac{dz}{zdx}=\frac{\alpha\beta}{x^{\beta+1}}
\]
and hence
\[
dz=\frac{a\alpha\beta x^{n-\beta-1}u^m}{v}dx.
\]

30. If we now write $k$ in place of $n-\beta-1$, so that $k$ and $m$ denote
any positive or negative numbers, since $n=k+\beta+1$ it will always be
\[
\frac{\int x^ku^m dx}{v}=\frac{1}{\alpha\beta}\cdot \frac{x^{k+\beta+1}u^m}{v}
\]
and hence if $\frac{x^ku^m}{v}=y$ and $y$ is seen as the ordinate of a curve,
its area will be
\[
\int ydx=\frac{1}{\alpha\beta} \cdot yx^{\beta+1},
\]
as long as $x$ is infinitely small, which is all the more noteworthy as
so far no way has been found for doing these integrations.

\end{document}